\documentclass{elsart}
\usepackage{amsfonts,amssymb,amsmath}

\DeclareMathOperator{\Aut}{Aut}
\DeclareMathOperator{\Gal}{Gal}
\DeclareMathOperator{\GL}{GL}
\DeclareMathOperator{\PGL}{PGL}
\DeclareMathOperator{\Br}{Br}
\DeclareMathOperator{\Frob}{Frob}
\newcommand{\kb}{\overline k}
\newcommand{\cM}{{\mathcal M}}
\newcommand{\Z}{{\mathbb Z}}
\newcommand{\bP}{{\mathbb P}}
\newcommand{\bA}{\bar A}
\newcommand{\ba}{\bar a}
\DeclareMathOperator{\charc}{char}

\begin{document}

\begin{frontmatter}

\title{Field of moduli and field of definition for curves of genus
  $2$}
\author{Gabriel Cardona\thanksref{BFM}}
\ead{gabriel.cardona@uib.es}
\address{Dept. Ci\`encies Matem\`atiques i Inf.,
  Universitat de les Illes Balears,\\
  Ed. Anselm Turmeda, Campus UIB.\\
  Carretera Valldemossa, km. 7.5,\\
  07071 -- Palma de Mallorca, Spain.
  }

\author{Jordi Quer\thanksref{BFM}}
\ead{quer@ma2.upc.es}
\address{Dept. Matem\`atica Aplicada II,
  Universitat Polit\`ecnica de Catalunya,\\
  Campus Nord, Ed. A0, Desp. 202.\\
  Jordi Girona, 1--3,\\
  08034 -- Barcelona, Spain}

\thanks[BFM]{Supported by BFM2000-0794-C02-02 and
  HPRN-CT-2000-00114}
\date{26 June 2002}

\begin{abstract}
  Let $\cM_2$ be the moduli space that classifies genus $2$ curves.
  If a curve $C$ is defined over a field $k$,
  the corresponding moduli point $P=[C]\in\cM_2$ is defined over $k$.
  In \cite{mestre91:_const}, Mestre solves the converse problem for
  curves with $\Aut(C)\simeq C_2$.
  Given a moduli point defined over $k$, Mestre finds an obstruction to
  the existence of a corresponding curve defined over $k$,
  that is an element in $\Br_2(k)$ not always trivial.
  In this paper we prove that for all the other possibilities of $\Aut(C)$,
  every moduli point defined over $k$
  is represented by a curve defined over $k$. We also give
  an explicit construction of such a curve in terms of the
  coordinates of the moduli point.
\end{abstract}

\begin{keyword}
  genus $2$ curves, moduli space, field of moduli, field of definition
\end{keyword}

\end{frontmatter}

\section{Preliminaries on genus $2$ curves}
\label{sec:prel-genus-2}

\paragraph*{Curves of genus $2$, sextic forms and invariants.}
\label{sec:curves-genus-2}

Let us fix $k$ a perfect field of characteristic $\charc k\neq 2$,
$\kb$ an algebraic closure, and
$G_k=\Gal(\kb/k)$ the absolute Galois group of $k$. If $C/k$ is a
genus $2$ curve, then it always
admits an affine model
\[y^2=f(x),\]
where $f(x)\in k[x]$ is a polynomial of degree $5$ or  $6$ without
multiple roots. We will call this model a hyperelliptic model for
$C$.

The classification of genus $2$ curves up to isomorphism was performed
by Clebsch and Bolza (cf. \cite{clebsch72:_theor_binaer_algeb_formen}
and \cite{bolza88}) using the classification of binary sextic forms up
to linear equivalence. Given a genus $2$ curve with hyperelliptic
model $y^2=f(x)$, the associated sextic form is given by
\[F(x_1,x_2)=x_2^6f(x_1/x_2).\]

The classification of binary sextic forms uses Hilbert's
theory of algebraic invariants
(cf. \cite{hilbert93:_theor_algeb_invar}).
An invariant is a polynomial expression $I\in k[a_0,\dots,a_6]$
in the coefficients of a generic sextic form
\[F(x_1,x_2)=a_0x_1^6+a_1x_1^5x_2+\dots+a_5x_1x_2^5+a_6x_2^6,\]
that
under a change of variables
\[(x_1,x_2)\mapsto(mx_1+nx_2,px_1+qx_2)\]
associated to
$M=\left(\begin{smallmatrix}m&n\\p&q\end{smallmatrix}\right)\in\GL_2(\kb)$
changes by a
fixed power of $\det M$.
Covariants are polynomial expressions in $k[a_0,\dots,a_6,x_1,x_2]$
with the same property.
The homogeneous degree in $x_i$ is called the order of the covariant,
while the degree in $a_i$ is called the degree of the covariant;
in particular, invariants are just covariants of order zero.
Given an invariant $I$ and a genus $2$ curve
$C$ or a sextic form $F$, we will denote by $I(C)$ or $I(F)$ the
value that takes the invariant when evaluated at this curve or form.

In terms of invariants, the classification goes as follows.
Two sextic forms $F$ and $F'$ are linearly equivalent if, and only if,
there exists $r\in \kb^*$ such that for every invariant $I$, one has
$I(F)=r^d I(F')$, where $d$ is the degree of $I$. As a consequence,
the fact that an invariant annihilates when evaluated to a certain
curve depends only on the isomorphism class of the curve.

The algebra generated by the invariants can be generated by five of them,
with degrees $2$, $4$, $6$, $10$ and $15$;
this implies that the condition of linear equivalence of forms
needs only to be tested for this generating set.
Among these generators, there is, up to equivalence, a single algebraic
relation which involves terms of degree $30$; namely, if $R$ is a generator
of degree $15$, then this relation states that $R^2$ is
a polynomial expression in the four generators of even degree.
As a consequence, the condition of equivalence
needs not to be tested for the invariant $R$.
Indeed, since the even degree invariants determine $R$ up to sign, after replacing $r$
with $-r$ if necessary, the condition holds for $R$ if it
holds for the even degree generators.

Therefore, one can reduce the
former condition on linear equivalence to test it
only for the four generator invariants of even degree.
In the literature appear several such sets of generators:
Clebsch invariants, that we will note as $c_2,c_4,c_6,c_{10}$,
and are noted $A,B,C,D$ in \cite{clebsch72:_theor_binaer_algeb_formen} and
\cite{mestre91:_const},
Igusa invariants, that we will note as $I_2,I_4,I_6,I_{10}$,
and are noted $A,B,C,D$ in \cite{igusa60:_arith_variet_modul_genus_two},
and the last set we will mention is that of Igusa arithmetic invariants
$J_2,J_4,J_6,J_{10}$
introduced by Igusa in \cite{igusa60:_arith_variet_modul_genus_two}.
All these generators have rational coefficients
when expressed as polynomials in the coefficients of a generic form.
However, while the invariants $J_i$ reduce well in any characteristic,
the $I_i$ fail in characteristics $2$ and $3$ and the $c_i$ in
characteristics $2$, $3$ and $5$.

The invariant $I_{10}=2^{12}J_{10}$ is just the discriminant of the form $F$.
The forms $F$ corresponding to genus 2 curves are therefore
those with $J_{10}(F)\neq0$.

\paragraph*{Absolute invariants and the moduli space.}
\label{sec:absol-invar-moduli}

Absolute invariants are quotients of invariants of equal degree.
The classification of genus $2$ curves up to isomorphism
can be stated in terms of absolute invariants:
Two curves $C$ and $C'$ are $\overline k$-isomorphic if, and only if,
for every absolute invariant $j$, one has $j(C)=j(C')$.
From this fact, Igusa constructed the
variety $\cM_2$ that classifies genus $2$ curves up to isomorphism as
a $3$-dimensional affine variety over $\Z$
(cf. \cite{igusa60:_arith_variet_modul_genus_two}).

In this moduli space, the classes of curves for which all invariants of
degree $2$ don't annihilate, which is equivalent to have $J_2(C)\neq 0$,
form an open set and we can take as coordinates for the moduli point
the absolute invariants
\[j'_1=\frac{J_2^5}{J_{10}},\qquad
j'_2=\frac{J_2^3J_4}{J_{10}},\qquad
j'_3=\frac{J_2^2J_6}{J_{10}}.\]
Similarly, for the set of curves for which all invariants of degree $2$
annihilate, but there exists a non-vanishing invariant of degree $4$,
that is $J_2=0$ and $J_4\neq 0$, one can take
\[
j''_1=0,\qquad
j''_2=\frac{J_4^5}{J_{10}^2},\qquad
j''_3=\frac{J_4J_6}{J_{10}}.\]
Finally, for the remaining points, one can take
\[j'''_1=0,\qquad
j'''_2=0,\qquad
j'''_3=\frac{J_6^5}{J_{10}^3}.\]

In must be noted that a point $P=[C]\in\cM_2$ is defined over a field $k$
if, and only if, its coordinates are in $k$, which amounts to say
that all absolute invariants, when particularized to $C$, are in
$k$. This way, the set of  points on $\cM_2(k)$ defined over $k$ is in
bijection with the set of $3$-tuples
$(j_1,j_2,j_3)$, $j_i\in k$, where
\[j_i=\begin{cases}
  j'_i, & \text{if $j_1\neq 0$,}\\
  j''_i, & \text{if $j_1=0$, $j_2\neq 0$,}\\
  j'''_i, & \text{if $j_1=j_2=0$.}
\end{cases}
\qquad i=1,2,3.
\]

At first sight, it could seem that to have
absolute invariants defined over $k$ does not guarantee that one can find
a set of corresponding generator invariants $J_2,J_4,J_6,J_{10},R\in k$
giving rise to them.
The following lemma proves the existence of such invariants.

\begin{lem}\label{thm:invariants-over-field}
  For every point $P\in \cM_2(k)$, there exist invariants defined over $k$ that
  give the absolute invariants which are coordinates for this point.
\end{lem}

\begin{pf}
  One can take
  \[(J_2,J_4,J_6,J_{10})=
  \begin{cases}
    (j_1,j_1j_2,j_1^2j_3,j_1^4),&\text{if $j_1\neq 0$,}\\
    (0,j_2,j_2j_3,j_2^2),&\text{if $j_1=0$, $j_2\neq 0$,}\\
    (0,0,j_3^2,j_3^3),&\text{if $j_1=j_2=0$,}
  \end{cases}
  \]
  and verify that with this choices, we get the values for the
  absolute invariants as defined above.

  We can also make a choice in such a way that $R$ is defined over $k$.
  Since $R^2$ can be written as a polynomial expression in the invariants $J_i$,
  we have that $R$ is of the form $R=\sqrt{d}R_0$ with $R_0\in k$.
  We can then take $r=\sqrt d$ and construct another set of invariats
  \begin{align*}
    (J'_2,J'_4,J'_6,J'_{10},R')&=
    (r^2 J_2,r^4 J_4,r^6 J_6,r^{10} J_{10},r^{15} R)\\
    &=
    (d J_2,d^2 J_4,d^3 J_6,d^5 J_{10},d^8 R_0),
  \end{align*}
  equivalent to the former and all of them defined over $k$.
  \qed
\end{pf}

\paragraph*{Group of automorphisms.}
\label{sec:moduli-space-group}

In \cite{bolza88}, Bolza gives the different
possibilities for the reduced group of automorphisms
$\Aut'(C)=\Aut(C)/\langle \imath\rangle$ of a genus $2$ curve,
where $\imath$ is the hyperelliptic involution.
In \cite{cardona99:_jacob_gl} the corresponding structures for
the full group $\Aut(C)$ are given.
If $\charc k\neq 2,3,5$, the posible automorphism groups of genus $2$ curves are
\[C_2,\ \ V_4,\ \ D_8,\ \ D_{12},\ \ 2D_{12},\ \ \tilde S_4,\ \
C_{10}.\]
The generic case corresponds to $\Aut(C)\simeq C_2$.
The curves having an involution different from the hyperelliptic involution form a
surface in $\cM_2$;
the curves with group of automorphisms
isomorphic to $D_8$ or $D_{12}$ describe two curves in this surface,
and they intersect in the single point corresponding to curves with
$\Aut(C)\simeq 2D_{12}$; the curve corresponding to $D_8$ has another
distinguished point corresponding to the case $\tilde S_4$. Finally,
the curves with group of automorphisms cyclic of order $10$ correspond
to an isolated point outside the surface of curves with non-hyperelliptic involutions.


The subvariety of curves admitting a non-hyperelliptic
involution is characterized by the condition $R=0$,
while the single moduli point of curves with group of automorphisms
isomorphic to $C_{10}$ is
characterized by the condition $J_2=J_4=J_6=0$, that is,
$(j_1,j_2,j_3)=(0,0,0)$.

In characteristic $3$, there is no curve with group of automorphisms isomorphic
to $2D_{12}$. In characteristic $5$, the three isolated points reduce
to a single one, for which the reduced group of automorphisms is
isomorphic to $\PGL(2,5)$.

\paragraph*{Field of moduli and field of definition.}
\label{sec:field-moduli}

Let $C$ be any curve defined over $\kb$.
We say that $k$ is a field of moduli for $C$
if $C\simeq {}^\sigma C$ for every $\sigma\in G_k$.
In the case of genus $2$ curves, this condition is equivalent to have the
absolute invariants defined over $k$ and, therefore,
to say that the point of moduli $P=[C]\in\cM_2$ is defined over $k$.

If $C$ is defined over $k$, then $P=[C]$ is defined over $k$.
The problem we want to solve is the converse:
Given a point $P\in\cM_2(k)$,
decide whether there is a curve representing $P$ defined over $k$.
In an equivalent formulation, we are interested in the following
Galois descent problem:
Given a genus $2$ curve defined over $\kb$ having $k$ as a field of moduli,
decide if there is an isomorphic curve defined over $k$.

The answer to this question is summarized in the following theorem,
the first part of which is proved by Mestre in \cite{mestre91:_const},

\begin{thm}\label{thm:teorema-gran}
  Let $P=[C]\in \cM_2(k)$.
  \begin{enumerate}
  \item (Mestre) If $\Aut(C)\simeq C_2$, there exists an obstruction to
    the existence of a curve $C'$ defined over $k$ and isomorphic to $C$,
    which is an element in $\Br_2(k)$ not necessarily trivial.
  \item If $\Aut(C)\not\simeq C_2$, there always exists $C'/k$ with $C\simeq C'$.
  \end{enumerate}
\end{thm}

The remaining sections of the paper contain the proof of this theorem,
split into several cases corresponding to different structures for $\Aut(C)$.
We remark that all proofs are constructive in the sense that they give an algorithm
for constructing a curve defined over $k$ starting from the coordinates of a moduli
point $P\in\cM_2(k)$. In the generic case $\Aut(C)\simeq C_2$, this is of course 
only possible when the obstruction is trivial, and for the construction one needs
to find a $k$-rational point on a conic defined over $k$; for all remaining cases
an expression for such a curve in terms of invariants is given.

We begin by reviewing the proof of Mestre for $\Aut(C)\simeq C_2$ since the way
we solve the case $\Aut(C)\simeq V_4$ is just a variation of his construction.

\section{The case $\Aut(C)\simeq C_2$}
\label{sec:case-autcisom-c_2}

Let $C:y^2=f(x)$ be a genus $2$ curve defined over $\kb$ with absolute
invariants defined over $k$ and $\Aut(C)\simeq C_2$.
We assume $\charc k\neq 2,3,5$.
Let $F(x_1,x_2)=x_2^6f(x_1/x_2)$ be the associated sextic.
Mestre considers a set of invariants and
covariants defined in Table \ref{old_invariants}.
We refer to \cite{mestre91:_const} for the definition and construction
of the covariants $(G,H)_k$ obtained from other covariants $G$ and $H$
by means of the so called \"uberschiebung (transvectant) operation.

\begin{table}
  \caption{Covariants suitable for the case $\Aut(C)\simeq C_2$
    \label{old_invariants}}
  \[
  \begin{array}{r@{\, =\, }l l l}\hline
    \multicolumn{2}{l}{\text{covariants}}&\text{order}&\text{degree}\\ \hline
    i & (F,F)_4 & 4 & 2 \\
    \Delta & (i, i)_2 & 4 & 4 \\
    Y_1 & (F,i)_4 & 2 & 3 \\
    Y_2 & (i,Y_1)_2 & 2 &5 \\
    Y_3 & (i,Y_2)_2 & 2 & 7 \\
    c_2 & (F,F)_6 & 0 & 2 \\
    c_4 & (i,i)_4 & 0 & 4 \\
    c_6 & (i,\Delta)_4 & 0 & 6 \\
    c_{10} & (Y_3,Y_1)_2 & 0 & 10 \\
    X_1 & (Y_2,Y_3)_1 & 2 & 12 \\
    X_2 & (Y_3,Y_1)_1 & 2 & 10 \\
    X_3 & (Y_1,Y_2)_1 & 2 & 8 \\
    A_{ij}& (Y_i,Y_j)_2 & 0 & 2(i+j+1) \\
    a_{ijk} & (F,Y_i)_2(F,Y_j)_2(F,Y_k)_2 & 0 & 2(i+j+k+2) \\
    R & -(Y_1,Y_2)_1(Y_2,Y_3)_1(Y_3,Y_1)_1 & 0 & 15 \\ \hline
  \end{array}
  \]
\end{table}

These covariants satisfy the relations
\begin{align*}
  \sum_{i,j=1}^3 A_{ij}X_iX_j&=0,\\
  \sum_{i,j,k=1}^3 a_{ijk}X_iX_jX_k&=R^3F,
\end{align*}
which should be understood as identities in $k[x_1,x_2]$.

Consider now the conic and cubic defined by the equations
\begin{align*}
L&:\sum_{i,j=1}^3 A_{ij}Y_iY_j,\\
M&:\sum_{i,j=1}^3a_{ijk}Y_iY_jY_k.
\end{align*}
Since the coefficients $A_{ij}$ and $a_{ijk}$ are invariants of even
degree, they can be given as a polynomial expression in the
invariants $J_i$. Moreover, after lemma
\ref{thm:invariants-over-field}, given a moduli point, one can find
$J_i$ over $k$ that represent the point and, therefore,
give equations for $L$ and $M$ with coefficients in $k$.
We remark that, in \cite{mestre91:_const}, Mestre makes linear changes
in these equations to put the coefficients
in terms of absulute invariants in order to have the varieties defined over $k$;
however, after lemma \ref{thm:invariants-over-field}, this is no longer necessary.

In this context, Mestre proves the following theorem.

\begin{thm}
  With the notations above, there exists a curve of genus $2$ defined
  over $k$ that represents the moduli point if, and only if,
  $L(k)\neq\emptyset$.
\end{thm}

Moreover, if the condition is satisfied, one can find a
model for $C$ defined over $k$ as follows. With a rational point $P\in
L(k)$ one constructs a $k$-isomorphism
\[(x_1,x_2)\mapsto(T_1(x_1,x_2),T_2(x_1,x_2),T_3(x_1,x_2)):\bP^1\to
L.\]
After replacing $Y_i$ with $T_i$ in the equation for the cubic $M$,
one obtains a sextic form with coefficients on $k$ and with the
desired absolute invariants.

\begin{rem}\label{remarca}
  This construction works well in every characteristic different from
  $2$, $3$ and $5$. This is because in the construction of the
  invariants $A_{i,j}$ and $a_{i,j,k}$ and the covariants $X_i$ one
  introduces denominators that are multiples of these
  primes. Nevertheless, if one multiplies this covariants and
  invariants by suitable integers, one gets another set of covariants
  and invariants with integer coefficients.

  With this modifications,
  the conic and cubic can be defined in every characteristic, and
  their coefficients can be expressed in terms of the invariants
  $J_i$. This way, one can make the method of construction work for
  characteristic $3$. However, the conic one obtains is always
  degenerated in characteristic $5$ and therefore this method
  fails to work in this characteristic.
\end{rem}

\section{The case $\Aut(C)\simeq V_4$}
\label{sec:case-autcisom-v_4}

Let $C:y^2=f(x)$ be a genus $2$ curve defined over $\kb$ with absolute
invariants defined over $k$ and $\Aut(C)\simeq V_4$.
As in the previous section we assume $\charc k\neq 2,3,5$,
although the remark \ref{remarca} also holds here.

The conic constructed by Mestre from
the absolute invariants of $C$ is in this case degenerated.
We mimic the construction of Mestre using
another pair of conic and cubic to solve the problem.

Let us introduce another set of invariants and covariants defined by
the expressions in Table \ref{new_invariants}.
\begin{table}
  \caption{Covariants suitable for the case $\Aut(C)\simeq V_4$
    \label{new_invariants}}
  \[
  \begin{array}{r@{\, =\, }l l l}\hline
    \multicolumn{2}{l}{\text{covariants}}&\text{order}&\text{degree}\\
    \hline
    \bar Y_1 & (F,i)_4 & 2 & 3 \\
    \bar Y_2 & (i,Y_1)_2 & 2 &5 \\
    \bar Y_3 & (Y_1,Y_2)_2 & 2 & 8 \\
    \bar X_1 & (\bar Y_2,\bar Y_3)_1 & 2 & 13 \\
    \bar X_2 & (\bar Y_3,\bar Y_1)_1 & 2 & 11 \\
    \bar X_3 & (\bar Y_1,\bar Y_2)_1 & 2 & 8 \\
    \bar A_{ij}& (\bar Y_i,\bar Y_j)_2 & 0 & 2(i+j+1)+\delta_3(i,j) \\
    \bar a_{ijk} & (F,\bar Y_i)_2(F,\bar Y_j)_2(F,\bar Y_k)_2 & 0 &
    2(i+j+k+2)+\delta_3(i,j,k) \\
    \bar R & -(\bar Y_1,\bar Y_2)_1(\bar Y_2,\bar Y_3)_1(\bar Y_3,\bar
    Y_1)_1 & 0 & 16 \\ \hline
  \end{array}
  \]
\end{table}
In this table,
$\delta_3$ stands for the number of its arguments equal to $3$.
Therefore, the invariants $\bar A_{13},\bar A_{23},\bar A_{31},\bar A_{32}$ have
odd degree, are multiple of $R$, and annihilate for curves
with non-hyperelliptic involutions.
Also, the invariants $\bar a_{ijk}$ with an odd number of subindexes equal to $3$
annihilate for these curves.
As for the rest, all of them have even degree and can be expressed as
polynomials in the invariants of even degree.
We get the following identities expressing the invariants $\bar A_{ij}$ and $\bar a_{ijk}$
in terms of Clebsch invariants and the $R$ invariant,
recalling that a permutation of indexes keeps the invariants unaltered:
\begin{align*}
  \bar A_{1,1}&=\frac{1}{3}(c_2 c_4 + 6{c_6}),\\
  \bar A_{1,2}&=\frac23(c_4^2 + c_2 c_6),\\
  \bar A_{1,3}&=0,\\
  \bar A_{2,1}&=\frac23(c_4^2 + c_2 c_6),\\
  \bar A_{2,2}&={c_{10}},\\
  \bar A_{2,3}&=0,\\
  \bar A_{3,1}&=0,\\
  \bar A_{3,2}&=0, \\
  \bar A_{3,3}&=\frac1{2\cdot3^2}(-4 c_4^4 - 8 c_2 c_4^2 c_6
   - 4 c_2^2 c_6^2+3 c_2 c_4 c_{10})
\end{align*}
and
\begin{align*}
\bar a_{1,1,1}&=\frac4{3^3\cdot5^2}(c_2^2 c_6 - 6 c_4 c_6 +
  9 c_{10}),\\
\bar a_{1,1,2}&=\frac2{3^3\cdot5^2}(2c_4^3 + 4 c_2 c_4 c_6 +
  12 c_6^2 + 3 c_2 c_{10}),\\
\bar a_{1,1,3}&=\frac{-1}{2\cdot3\cdot5^2}R,\\
\bar a_{1,2,2}&=\frac2{3^4\cdot5^2}(3 c_2 c_4^3 +
  4 c_2^2 c_4 c_6 +
  12 c_4^2 c_6 + 18 c_2 c_6^2 +
  9 c_4 c_{10}),\\
\bar a_{1,2,3}&=\frac{-1}{2^2\cdot3^2\cdot5^2}c_2R,\\
\bar a_{1,3,3}&=\frac1{3^5\cdot5^2}(-3c_2^2c_4^4 + 24c_4^5 -
  4c_2^3c_4^2c_6 +
  42c_2c_4^3c_6 + 6c_2^2c_4c_6^2 + 72c_4^2c_6^2 \\ &\qquad+
  36c_2c_6^3 + 27c_2c_4^2c_10 + 18c_2^2c_6c_{10} +
  54c_4c_6c_{10} - 162c_{10}^2),\\
\bar a_{2,2,2}&=\frac2{3^4\cdot5^2}(9 c_4^4 + 18c_2c_4^2c_6 +
  8 c_2^2 c_6^2+6c_4c_6^2 -
  9c_6c_{10}),\\
\bar a_{2,2,3}&=\frac{-1}{2^2\cdot3^2\cdot5^2}c_4R,\\
\bar a_{2,3,3}&=\frac1{3^5\cdot5^2}(3c_2c_4^5 + 10c_2^2c_4^3c_6 -
  6c_4^4c_6 +
  8c_2^3c_4c_6^2 + 6c_2c_4^2c_6^2 +
  24c_2^2c_6^3 \\ &\qquad    - 36c_4c_6^3 + 18c_4^2c_{10} +
  9c_2c_4c_6c_{10} -
  54c_6^2c_{10} -
  27c_2c_{10}^2),\\
\bar a_{3,3,3}&=\frac{1}{2^3\cdot3^3\cdot5^2}(-3c_2c_4^2-
4c_2^2c_6+6c_4c_6+18c_{10})R.
\end{align*}
For the set of invariants and covariants defined, one has identities
analogous to those used by Mestre to define $L$ and $M$. Namely, one
has
\[\sum_{i,j=1}^3 \bar A_{ij}\bar X_i\bar X_j=0\]
and
\[\sum_{i,j,k=1}^3 \bar a_{ijk}\bar X_i\bar X_j\bar X_k=\bar R^3F.\]

Let us consider the conic and cubic defined by
\begin{align*}
  \bar L&:\sum_{i,j}\bar A_{i,j} Y_iY_j,\\
  \bar M&:\sum_{i,j,k}\bar a_{i,j,k} Y_iY_jY_k.
\end{align*}
As in the case $\Aut(C)\simeq C_2$, one can construct $\bar L$ and
$\bar M$ defined over $k$ from a moduli point defined over $k$.

\begin{lem}
  With the notations as above, and considering a moduli point
  corresponding to group of automorphisms $V_4$,
  the conic $\bar L$ is non-degenerated.
\end{lem}

\begin{pf}
  The discriminant of $\bar L$ is an invariant; therefore, the fact
  that it does annihilate or not depends only on the moduli point.
  We can then suppose that we have a representative for the
  point with hyperelliptic model
  \[y^2=x^6+\alpha x^4+\beta x^2+1,\]
  with $\alpha,\beta\in\kb$. The fact that a genus $2$ curve with
  non-hyperelliptic involutions always is isomorphic to some curve of
  this type was proved by Bolza (cf. \cite{bolza88}).

  Computing explicitly the discriminant of $\bar L$ constructed from
  this model over $\kb$ we get a non-zero multiple of
  \[( \alpha - \beta )^4 ( \alpha^2 + \alpha\beta + \beta^2)^4
  ( -1125 + 4\alpha^3 + 110\alpha\beta - \alpha^2\beta^2 + 4\beta^3 )^4.\]
  It is now easy to check that if any of these terms annihilates, then
  $\Aut(C)$ contains $D_8$ or $D_{12}$.
  \qed
\end{pf}

\begin{thm}
  Let $C$ be a genus $2$ curve with group of automorphisms isomorphic
  to $V_4$ and field of moduli $k$. Then, $C$ is isomorphic to the
  curve defined over $k$ given by the model
  \begin{multline*}
    y^2=\Big(
    - \bA_{3,3} \ba_{1,1,1} P_1^3(x)
    -3 \bA_{3,3} \ba_{1,1,2} P_1^2(x) P_2(x) \\
    -3 \bA_{3,3} \ba_{1,2,2} P_1(x) P_2^2(x)
    +3 \bA_{2,2} \ba_{1,3,3}P_1 P_3^2(x)\\
    -\bA_{3,3}\ba_{2,2,2}P_2^3(x)
    +3 \bA_{2,2} \ba_{2,3,3}P_2(x)P_3^2(x)
    \Big),
  \end{multline*}
  where
  \begin{align*}
    P_1(x)&=-2\bA_{1,2}-2 \bA_{2,2}x, \\
    P_2(x)&=\bA_{1,1}-\bA_{2,2}x^2,\\
    P_3(x)&=\bA_{1,1}+2\bA_{1,2}x+\bA_{2,2}x^2;
  \end{align*}
  whose coefficients can be explicitly obtained from the invariants of $C$.
\end{thm}

\begin{pf}
  Using the previous construction, one gets a non-degenerated conic
  with coefficients in $k$.
  The conic $\bar L$ has a point defined on a quadratic extension of
  $k$; namely,
  \[P=\left(0,\sqrt{-\bar A_{3,3}},\sqrt{\bar A_{2,2}}\right).\]
  Using this point, one obtains the following parametrization for the conic
  \[
  \lambda\mapsto \left(
    P_1(\lambda)\sqrt{-\bA_{3,3}},P_2(\lambda)\sqrt{-\bA_{3,3}},
    P_3(\lambda)\sqrt{\bA_{2,2}}\right),
  \]
  where $P_i$ are the polynomials given in the statement of the theorem.
  Using the equation for the cubic $\bar M$ with $R=0$, since we are
  assuming $\Aut(C)\simeq V_4$, and the parametrization
  obtained, one obtains that $C$ is isomorphic to the curve given by
  $y^2=f(x)$, where $f$ is
\begin{multline*}
  \sqrt{-\bA_{3,3}}\Big(
  - \bA_{3,3} \ba_{1,1,1} P_1^3(x)
  -3 \bA_{3,3} \ba_{1,1,2} P_1^2(x) P_2(x) \\
  -3 \bA_{3,3} \ba_{1,2,2} P_1(x) P_2^2(x)
  +3 \bA_{2,2} \ba_{1,3,3}P_1(x) P_3^2(x)\\
  -\bA_{3,3}\ba_{2,2,2}P_2^3(x)
  +3 \bA_{2,2} \ba_{2,3,3}P_2(x)P_3^2(x)
  \Big).
\end{multline*}
Therefore, and up to isomorphism, a model for $C$ defined over $k$ is
given by the equation in the theorem.
\qed
\end{pf}

\begin{rem}
  The pair $\bar L$ and $\bar M$ can also be used as alternatives to
  Mestre's $L$ and $M$ in the case $\Aut(C)\simeq C_2$ for all moduli points outside the
  surface detemined by $\bar R=0$.
  However, when $R\neq 0$ the previous formula is not valid,
  and the construction needs to start with a point $P\in \bar L(k)$.

  If we construct the conics $L$ and $\bar L$ from invariants $c_i\in k$
  with, also, $R\in k$, then of course the obstructions, as elements
  in $\Br_2(k)$,
  to find a $k$-rational point on $L$ and on $\bar L$ coincide.
\end{rem}

\section{The case $\Aut(C)\simeq D_8,D_{12}$}
\label{sec:case-autcisom-d_8-d_12}

These cases have been studied in detail in \cite{cardonaquer}.
Therefore, we only state the results and refer to the mentioned paper for the proof,
except for the case of $\Aut(C)\simeq D_{12}$ in characteristic $3$,
not treated in \cite{cardonaquer}, for which we give the proof here.

\begin{thm}\label{thm:d8-sobre-kbarra}
  Let $C$ be a genus $2$ curve with $\Aut(C)\simeq D_8$ and field of
  moduli $k$. Then, $C$ is $\kb$-isomorphic to the curve defined over $k$
  \[y^2=x^5+x^3+tx,\]
  where $t\in k\setminus\{0,1/4\}$ is the absolute invariant of $C$ given
  by
  \[t=t(C)=\begin{cases}
     1+\dfrac{J_4}{J_2^2},& \text{if $\charc k=5$,}\\[10pt]
     \dfrac{-J_2^2}{J_4},& \text{if $\charc k=3$,}\\[10pt]
    \dfrac{8 c_6(6 c_4-c_2^2)+9c_{10}}{900 c_{10}},&\text{if $\charc
      k\neq 3,5$.}
  \end{cases}
  \]
\end{thm}

\begin{thm}\label{thm:d12-sobre-kbarra}
  Let $C$ be a genus $2$ curve with $\Aut(C)\simeq D_{12}$ and field of
  moduli $k$. If $\charc k\neq 3$, then $C$ is $\kb$-isomorphic to the
  curve defined over $k$
  \[y^2=x^6+x^3+t,\]
  where $t\in k\setminus\{0,1/4\}$ is the absolute invariant of $C$ given
  by
  \[t=t(C)=\begin{cases}
    -1-\dfrac{J_4}{J_2^2},& \text{if $\charc k=5$,}\\[10pt]
    \dfrac{3c_4c_6-c_{10}}{50c_{10}},&\text{if $\charc
      k\neq 3,5$.}
  \end{cases}\]
  If $\charc k=3$, then $C$ is $\kb$-isomorphic to the
  curve defined over $k$
  \[y^2=\frac{1}{t_*}x^6+x^4+x^2+1,\]
  where $t_*\in k^*$ is such that $t_*^3=t$ and $t\in k^*$ is the
  absolute invariant given by
  \[t=t(C)=\frac{-J_2^3}{J_6}.\]
\end{thm}

Consider this last $\charc k=3$ case.
In the same way as the proofs in \cite{cardonaquer},
we see that the group $\Aut(C)\subset \GL_2(\kb)$ isomorphic to
$D_{12}$ is $\GL_2(\kb)$-conjugated to the group $\langle U,V\rangle$ generated by
\[U=\begin{pmatrix}{-1}&1\\{-1}&0\end{pmatrix},\qquad
V=\begin{pmatrix}0&1\\1&0\end{pmatrix}.\]
This implies that $C$ is $\kb$-isomorphic to the curve
\[y^2=a_{6}x^{6}+a_{4}x^{4}+(a_{4}+2a_{6})x^{3}+a_{4}x^{2}+a_{6},\]
where $a_4,a_6\in\kb^*$, because otherwise the polynomial has multiple roots.
After changing this curve with the isomorphism associated to the matrix
$\sqrt{a_4}\left(\begin{smallmatrix}1&1\\0&1\end{smallmatrix}\right)$,
one obtains the isomorphic curve given by the equation
\[y^2=\frac{a_6}{a_4}x^6+x^4+x^2+1,\]
and $C$ is $\kb$-isomorphic to a the curve with hyperelliptic model
\[C_{t_*}:y^2=\frac{1}{t_*}x^6+x^4+x^2+1,\qquad t_*=\frac{a_6}{a_4}in\kb^*.\]
Computing the absolute invariant of the theorem one obtains that $t(C_{t_*})=t_*^3$.
Therefore, since we are in characteristic 3, $\Frob_3$ is a bijection,
and it follows that $C_{t_*}$ and $C_{t'_*}$ are isomorphic if, and only if, $t_*=t'_*$.
Moreover, if $C$ has $k$ as a field of moduli, then $t=t(C)$ belongs to $k$;
since $\Frob_3$ is an automorphism of $k$,
we can find $t_*\in k$ with $t=t_*^3$ and $C$ is isomorphic to the
curve defined over $k$ with equation
\[y^2=\frac{1}{t_*}x^6+x^4+x^2+1.\]

\section{The case $\Aut(C)\simeq 2D_{12},\tilde S_4,C_{10}$}
\label{sec:case-autcisom-2d_12-s4t-c10}

The remaining trivial cases correspond to three single points in $\cM_2$.
For the sake of completeness, we give below an explicit equation defined over $k$
for every moduli point.

\begin{thm}
  Let $C$ be a genus $2$ curve with group of automorphisms isomorphic to
  $2D_{12}$, $\tilde S_4$ or $C_{10}$. Then, $C$ is
  $\kb$-isomorphic to the curve given by
  \begin{align*}
    y^2&=x^6-1,& &\text{if $\Aut(C)\simeq 2D_{12}$,}\\
    y^2&=x^5-x,& &\text{if $\Aut(C)\simeq \tilde S_4$,}\\
    y^2&=x^5-1,& &\text{if $\Aut(C)\simeq C_2\times C_5$.}\\
  \end{align*}
\end{thm}


\begin{thebibliography}{1}

\bibitem{bolza88}
Oskar Bolza.
\newblock On binary sextics with linear transformations into themselves.
\newblock {\em Amer. J. Math}, 10:47--70, 1888.

\bibitem{cardona99:_jacob_gl}
Gabriel Cardona, Josep Gonz{\'a}lez, Joan~Carles Lario, and Anna Rio.
\newblock On curves of genus $2$ with {J}acobian of ${{\rm {G}{L}}}\sb 2$-type.
\newblock {\em Manuscripta Math.}, 98(1):37--54, 1999.

\bibitem{cardonaquer}
Gabriel Cardona and Jordi Quer.
\newblock Curves of genus $2$ with group of automorphisms isomorphic to ${D_8}$
  or ${D_{12}}$.
\newblock Submitted, February 2002.

\bibitem{clebsch72:_theor_binaer_algeb_formen}
Alfred Clebsch.
\newblock {\em Theorie der Bin{\"a}rien Algebraischen Formen}.
\newblock Verlag von B. G. Teubner, Leipzig, 1872.

\bibitem{hilbert93:_theor_algeb_invar}
David Hilbert.
\newblock {\em Theory of {A}lgebraic {I}nvariants}.
\newblock Cambridge University Press, 1993.

\bibitem{igusa60:_arith_variet_modul_genus_two}
Jun-Ichi Igusa.
\newblock Arithmetic variety of moduli for genus two.
\newblock {\em Annals of Mathematics}, 72(3):612--649, 1960.

\bibitem{mestre91:_const}
Jean-Fran{\c{c}}ois Mestre.
\newblock Construction de courbes de genre $2$ \`a partir de leurs modules.
\newblock In {\em Effective methods in algebraic geometry (Castiglioncello,
  1990)}, pages 313--334. Birkh\"auser Boston, Boston, MA, 1991.

\end{thebibliography}
\end{document}